\documentclass[10.5pt]{article} 
\usepackage{fouriernc}
\usepackage{extarrows}
\usepackage{array}

\renewcommand{\arraystretch}{1.3}

\usepackage{amsfonts,amsmath,amssymb,amsthm,amsxtra,graphicx,latexsym}
\numberwithin{equation}{section}

\def\sct{\section}
\def\ssct{\subsection}

\def\ba{\sla\begin{array}}\def\ea{\end{array}}     
\def\be{\begin{equation}}\def\ee{\end{equation}}   
\def\bgd{\begin{aligned}}\def\egd{\end{aligned}}   
    \def\lgd{\lt\{\bgd}\def\rgd{\egd\rt.}          
\def\bgt{\begin{gathered}}\def\egt{\end{gathered}} 

\def\bcc{\begin{center}}\def\ecc{\end{center}}
\def\bmc{\begin{multicols}}\def\emc{\end{multicols}}
\def\bmpg{\begin{minipage}}\def\empg{\end{minipage}}
\def\btb{\begin{tabular}}\def\etb{\end{tabular}}
\def\bfl{\begin{flushleft}}
\def\efl{\end{flushleft}}

\def\lb{\label}
\def\rf{\ref}

\def\<{\langle}
\def\>{\rangle}

\def\bs{\setminus}
\def\cd{\cdots}

\def\dtp{\dotplus}

\def\ev{\equiv}

\def\fa{\forall}
\def\gsl{\geqslant}

\def\ift{\infty}

\def\lsl{\leqslant}
\def\lls{\lim\limits}

\def\lt{\left}\def\rt{\right}

\def\op{\oplus}

\def\pt{\partial}
\def\ra{{\rightarrow}}

\def\ts{\times}
\def\udb{\underbrace}

\def\wtd{\widetilde}
\def\xeq{\xlongequal}

\def\aa{\alpha}
\def\bb{\beta}
\def\ga{\gamma}
\def\dl{\delta}

\def\tht{\theta}

\def\lm{\lambda}

\def\vf{\varphi}

\def\N{\mathbb N}

\def\R{\mathbb R}
\def\Z{\mathbb Z}

\def\cE{\mathcal E}

\def\cH{\mathcal H}

\def\cP{\mathcal P}
\def\rd{\textrm d}

\makeatletter
\def\ExtendSymbol#1#2#3#4#5{\ext@arrow 0099{\arrowfill@#1#2#3}{#4}{#5}}
\def\RightExtendSymbol#1#2#3#4#5{\ext@arrow 0359{\arrowfill@#1#2#3}{#4}{#5}}
\def\LeftExtendSymbol#1#2#3#4#5{\ext@arrow 6095{\arrowfill@#1#2#3}{#4}{#5}}
\makeatother

\newcommand\mRa[2][]{\ExtendSymbol{~}{=}{\Rightarrow}{#1}{#2}}

\DeclareSymbolFont{EUEX}{U}{euex}{m}{n}
\DeclareSymbolFont{euexlargesymbols}{U}{euex}{m}{n}
\DeclareMathSymbol{\intop}{\mathop}{euexlargesymbols}{"52}
\def\int{\intop\nolimits}

\DeclareSymbolFont{mathdesign}{OMX}{mdbch}{m}{n}
\DeclareMathSymbol{\intmd}{\mathop}{mathdesign}{"52}
\DeclareMathSymbol{\ointmd}{\mathop}{mathdesign}{"49}
\def\int{\intop\nolimits}

\DeclareSymbolFont{ugmL}{OMX}{mdugm}{m}{n}
\SetSymbolFont{ugmL}{bold}{OMX}{mdugm}{b}{n}
\DeclareMathAccent{\wideparen}{\mathord}{ugmL}{"F3}

\def\QEDclosed{\mbox{\rule[0pt]{0.7ex}{1.3ex}}} 

\def\QED{\QEDclosed} 

\def\hs{\hspace}    

\def\nid{\noindent} 
\def\nn{\nonumber}  
\def\ras{\renewcommand{\arraystretch}}
\def\sla{\setlength{\arraycolsep}{2pt}} 
\def\sss{\scriptscriptstyle}

\def\text{\textrm}
\def\vs{\vspace}

\def\ncP{^\#\cP}

\def\bbm{\begin{bmatrix}}\def\ebm{\end{bmatrix}}
\def\bsm{\begin{smallmatrix}}\def\esm{\end{smallmatrix}}

\def\sI{\text{I}}
\def\sII{\text{I\!I}}

\newtheoremstyle{plain}{0.6ex}{0.6ex}{}{}{}{}{0.5em}{\bf\thmname{#1}~\thmnumber{#2}.~\thmnote{(#3)}}
\newtheoremstyle{definition}{0.6ex}{0.6ex}{}{}{}{}{0.5em}{\bf\thmname{#1}~\thmnumber{#2}.~\thmnote{(#3)}}

\theoremstyle{plain}
\newtheorem{que}{Question}
\newtheorem{thm}{Theorem}
\newtheorem{cor}{Corollary}

\newtheorem{prop}{Proposition}[section]
\newtheorem{lem}{Lemma}
\newtheorem{clm}{Claim}

\theoremstyle{definition}
\newtheorem{Def}{Definition}
\newtheorem{exa}{Example}
\newtheorem{rem}{Remark}

\def\bQ{\begin{que}}\def\eQ{\end{que}}
\def\bT{\begin{thm}}\def\eT{\end{thm}}
\def\bcr{\begin{cor}}\def\ecr{\end{cor}}

\def\bP{\begin{prop}}\def\eP{\end{prop}}
\def\blm{\begin{lem}}\def\elm{\end{lem}}
\def\bcl{\begin{clm}}\def\ecl{\end{clm}}

\def\bD{\begin{Def}}\def\eD{\end{Def}}
\def\bE{\begin{exa}}\def\eE{\end{exa}}
\def\bR{\begin{rem}}\def\eR{\end{rem}}

\def\pf{\noindent{\it Proof.~~}}

\def\rT{Theorem~\rf}

\allowdisplaybreaks

\begin{document}
\title{
\vspace{0.5in} {\bf\Large Rotational solutions of prescribed period of Hamiltonian systems on $\R^{2n-k}\ts T^k$}} 
\author{{\bf\large Hui Qiao\vspace{1mm}}\\
{\it\small School of Mathematics and Statistics}\\
{\it\small Wuhai University, Wuhan 430072},\\
{\it\small The People's Republic of China}\\
{\it\small e-mail: qiaohuimath@mail.whu.edu.cn}}\vspace{2mm}
\date{
}
\maketitle
\begin{center}
{\bf\small Abstract}

\vspace{3mm} \hspace{.05in}\parbox{4.5in}
{\small In this paper, we consider Hamiltonian systems on $\R^{2n-k}\ts T^k$.
Multiple rotational solutions are obtained.}
\end{center}

\sct{Introduction and main results}

In this paper we study the multiplicity of rotational solutions for
Hamiltonian systems
\be\dot z=JH'(z),\quad J=\lt(\ba{cc} 0 & -I_n \\ I_n & 0\ea\rt).\lb{eq:HS}\tag{HS}\ee

For $1\lsl k\lsl 2n-1$, let
\be z=(z_{\sI},z_{\sII}),~z_{\sI}=(z_1,\cd,z_{2n-k}),
~z_{\sII}=(z_{2n-k+1},\cd,z_{2n}),\ee
and
\be H_{z_{\sI}}=(H_{z_1},\cd,H_{z_{2n-k}}),~H_{z_{\sII}}=(H_{z_{2n-k+1}},\cd,z_{2n}).\lb{eq:HzI,HzII}\ee
We make the following basic hypothesis on the Hamiltonian

(H0)\quad $H\in C^1(\R^{2n},\R)$.

(H1)\quad $H(z_{\sI},z_{\sII}+v)=H(z_{\sI},z_{\sII}),~\fa z=(z_{\sI},z_{\sII})\in\R^{2n},~v\in\Z^k$.

(H2)\quad There exist constants $r>0,\,\mu>1$ such that
\[0<\mu H(z)\lsl z_{\sI}\cdot H_{z_{\sI}}(z),~\fa z=(z_1,\cd,z_{2n})\in\R^{2n},
~|z_{\sI}|^2\gsl r^2.\]

\nid A periodic solution $z(t)$ of \eqref{eq:HS} on $\R^{2n-k}\ts T^k$ satisfies
\be z(T)-z(0)=(0,v)\lb{eq:bdc}\ee
for period $T>0$ and vector $v\in\Z^k$. If $v=0$, solutions are contractible.
If $v\ne0$, it is called {\it rotational vector} and corresponding solutions are
called {\it rotational solutions}, which are non-contractible.

If $z(t)$ is a rotational solution of \eqref{eq:HS}, then $z(mt)$ is also such a
solution with rotational vector $mv$. We restrict choices of rotational vectors to
$\Z^k_1$ defined as follows.

\bD A rotational vector $v=(v_1,v_2,\cd,v_k)\in\Z^k\bs\{0\}$ is called {\it prime}, if
one of its coordinates equals to $1$ while others are zero, or two of them are relatively
prime. Denotes by
\[\Z^k_1=\text{the set of all prime rotational vectors in}~\Z^k.\]\eD

In this paper, we consider the following boundary value problem
\be\lt\{\bgd\dot z(t)&=JH'(z(t)),\\ z(T)&=z(0)+(0,v),\egd\rt.\lb{eq:fpH.z}\ee

For given $T>0$ and $v\in\Z^k_1$. Denote by $\cP_H(T,v)$ the set of distinct solutions.
Main results in this paper are the following theorems.

\bT\lb{thm:rot.fp.k=n}\it Assume that $k=n$ and $H$ satisfies (H0)\,--\,(H2). For every $T>0$ and $v\in\Z^k_1$, we have
\[\ncP_H(T,v)\gsl k.\]
\eT

\bT\lb{thm:rot.fp.k>n}\it Assume that $k>n$ and $H$ satisfies (H0)\,--\,(H2) and


(H3)\quad There exist positive numbers $a,b$ and $s<\mu-\frac{1}{2}$ such that
\[|H_{z_{\sII}}(z)|\lsl a|z_{\sss\text{I}}|^s+b,~~\fa z=(z_{\sI},z_{\sII})\in\R^{2n-k}\ts\R^k.\]

\nid For every $v\in\Z^k_1\cap\lt(\{0\}\ts\Z^{2n-k}\ts\{0\}\rt)$ with $0\in\R^{k-n}$ and $T\in I(s)$, we have
\[\ncP_H(T,v)\gsl k,\]
where the interval $I(s)$ is defined as follows
\[I(s)=\lgd (0,+\ift),&\quad s<\frac{\mu}{2},\\ (0,\dl),&\quad \frac{\mu}{2}\lsl s<\mu-\frac{1}{2},
\rgd\]
and $\dl$ is sufficiently small.
\eT






\sct{Variational settings}

Make substitution
\be z(Tt)=x(t)+t(0,v),\lb{eq:z=x+v}\ee
then solutions of \eqref{eq:fpH.z} are in one to one correspondence with $1$-periodic
solutions of the following problem:
\be\lt\{\bgd\dot x(t)+(0,v)&=TJH'(x(t)+t(0,v)),\\ x(1)&=x(0).\egd\rt.\lb{eq:fpH.x}\ee


We choose the Hilbert space $E=W^{\frac{1}{2},2}(S^1,\R^{2n})$ with its inner product and norm defined by
\be\<x,y\>=x_0\cdot y_0+\sum_{j\in\Z}2\pi|j|x_j\cdot y_j,~\fa x,y\in E,\ee
and
\be ||x||=\<x,x\>^{\frac{1}{2}},~\lb{eq:defnjnorm}\ee
where
\[x(t)=\sum_{j\in\Z}e^{2\pi jtJ}\xi_j,~\xi_j\in\R^{2n},\]
is the Fourier series expansion of $x$. For $x\in E\cap C^\ift$ we define
\be A(x)=-\frac{1}{2}\int^1_0J\dot x(t)\cdot x(t)\rd t.\lb{eq:defA.1}\ee
It can be extended to the whole space $E$. Note that the space $E$ can be orthogonally decomposed as
\be E=E^+\op E^-\op E^0,\lb{eq:E.dec}\ee
where
\begin{align}
E^{\pm}&=\lt\{\sum_{j\in\Z\bs\{0\}}e^{2\pi jtJ}\xi_j~\bigg|~\pm j>0\rt\}\nn\\
&=\text{span}_\R\lt\{\ras{0.8}\ba{c}\sin(2\pi jt)e_i\mp\cos(2\pi jt)e_{i+n},\\
\cos(2\pi jt)e_i\pm\sin(2\pi jt)e_{i+n}\ea~\Big|~1\lsl i\lsl n,~j\in\N\rt\},\lb{eq:defE+-}\\
E^0&=\text{span}_\R\big\{e_1,\cd,e_{2n}\big\}=\R^{2n}.\lb{eq:defE0}\end{align}
Denote by $P^0$ and $P^\pm$ the orthogonal projections from $E$ onto $E^0$ and $E^\pm$
respectively. Let
\be x^0=P^0x,~x^\pm=P^\pm x.\lb{eq:defx0+-}\ee
Define a linear, bounded and selfadjoint operator
\be L=P^+-P^-:E\ra E.\lb{eq:defL}\ee
Then the functional $A:~E\cap C^\ift\ra\R$ defined by \eqref{eq:defA.1} can be
extended to $E$ as follows
\be A(x)=\frac{1}{2}\<Lx,x\>=\frac{1}{2}\lt(||x^+||^2-||x^-||^2\rt).\lb{eq:defA.2}\ee

Consider a Hamiltonian $H$ satisfying (H0),
(H1) and the growth condition

\be|H(z)|\lsl a|z_{\sI}|^s+b,~\fa z=(z_{\sI},z_{\sII}),~z_{\sI}\in\R^{2n-k},~z_{\sII}\in\R^k,\lb{eq:grtcd}\ee

\nid for some constants $a,b,s>0$. Then define a functional $B: E\ra\R$ by
\be B(x)=\int^1_0 \big[H(x(t)+t(0,v))+x(t)\cdot J(0,v)\big]\rd t.\lb{eq:defB}\ee
This functional is well defined, of class $C^1$ and its derivative $B'$ is compact.
Now we can define on $E$ the functional
\be \Phi(x)=A(x)-B(x).\lb{eq:defPhi}\ee

\bP $x\in E$ is a critical point of $\Phi$ iff $x$ is a solution of \eqref{eq:fpH.x}.\eP

For proofs of main theorems, we need different decompositions instead of \eqref{eq:E.dec}. Note that $E^0$
can be decomposed as
$E^0=E^0_{\sI}\op E^0_{\sII}$, where
\begin{align} &E^0_{\sI}=\text{span}\big\{e_1,\cd,e_{2n-k}\big\},\\
&E^0_{\sII}=\text{span}\big\{e_{2n-k+1},\cd,e_{2n}\big\}.\lb{eq:defE0}\end{align}

We define subspaces $X,Y,\cH$ of $E$ such that

\be \cH=X\dtp Y=E^+\op E^-\op E^0_{\sI}.\lb{eq:defXYZcH}\ee
Then the space $E$ can be decomposed as
\be E=\cH\op E^0_{\sII}.\lb{eq:E.dec'}\ee
Let
\begin{align}
&E_{\sI}=\lt\{(x_1,\cd,x_{2n-k},0)\in E\,\bigg|\,\int^1_0 x_i(t)\rd t=0\rt\},\\
&E_{\sII}=\lt\{(0,x_{2n-k+1},\cd,x_{2n})\in E\,\bigg|\,\int^1_0 x_i(t)\rd t=0\rt\}.\end{align}
Subspace $X,Y$ can be defined according to various values of $k$ as follows.\\

{\bf Case 1.} $k=n$.
\be X=E^-\op E^0_{\sI},~Y=E_{\sII}.\lb{eq:XY:k=n}\ee

{\bf Case 2.} $k>n$.\quad
$x=(\udb{\ast,\cd,\ast}_{2n-k},\udb{\ast,\cd,\ast}_{k-n},\udb{\ast,\cd,\ast}_{2n-k},\udb{\ast,\cd,\ast}_{k-n})$.

\[X_1=\lt\{(p_1,\cd,p_{2n-k},0,q_1,\cd,q_{2n-k},0)\in E\,\big|\,\int^1_0 p_i(t)\rd t=\int^1_0 q_i(t)\rd t=0\rt\}.\]

\be X=\lt(E^-\cap X_1\rt)\op E^0_{\sI},~Y=E_{\sII}.\lb{eq:XY:k>n}\ee






We consider symmetries involved in the problem.

Firstly, since $H$ satisfies (H1), the functional $\Phi$ defined by \eqref{eq:defPhi} is $\Z^k$-invariant and
can be defined on
\be\cE=E/\Z^k=\cH\ts T^k.\ee

Secondly, P. Felmer \cite{pF92TMA} points out that the following $S^1$-action is free,
and the functional $\Phi$ is $S^1$-invariant.
\be(\tht\cdot x)(t)=x(t+\tht)+(0,\tht v),~x\in\cE.\lb{eq:s1act}\ee
where $\tht\in S^1=[0,1]/\{0,1\}$.



P. Felmer \cite{pF92JDE} essentially proves the following $S^1$-equivariant saddle point type theorem.

\bP[Theorem 1.1 of \cite{pF92JDE}]\lb{prop:sptt.fp}\it Assume that $\cE$ can be splitted as
$\cE=(X\dtp Y)\ts T^k$ with $X\neq\{0\}$. Let $I\in C^1(\cE,\R)$ be an
$S^1$-invariant functional of the form
\be I(z)=\frac{1}{2}\<Lz,z\>+B(z),\lb{eq:fp.defI}\ee
where
\par(I1)\quad $L:E\ra E$ is a linear, bounded and selfadjoint operator; $X$ is an invariant
subspace.
\par(I2)\quad $b\in C^1(\cE,\R)$ and $b'_z:\cE\ra\cH$ is compact.
\par(I3)\quad $I$ satisfies (PS) property.
\par(I4)\quad There exist constants $\aa<\bb$ and $\ga$ such that
\[I|_{\pt Q\ts T^k}\lsl\aa,~I|_{Y\ts T^k}\gsl\bb,~I|_{Q\ts T^k}\lsl\ga,\]
where
\be Q=\{x\in X~|~||x||\lsl R\}~~\text{and}~~\pt Q=\{x\in X~|~||x||=R\}.\lb{eq:defQ}\ee
Then $I$ possesses at least $k$ distinct critical points with critical values less than or equal to $\ga$.\eP

\sct{Proofs of main results}

Since $H$ satisfies (H0)-(H2), then
\be H(z)\gsl a_1|z_{\sI}|^\mu-a_2,~\fa z=(z_{\sI},z_{\sII})\in\R^{2n-k}\ts\R^k,\ee
where
\be a_1=\min_{|z_\sI|\gsl r}\frac{H(z)}{|z_\sI|^\mu},\ \ a_2=\max_{|z_\sI|\lsl r}|H(z)|.\ee
Choose undetermined constants $K_2>K_1\gsl r$ and a function $\chi\in C^\ift(\R^+,\R^+)$ such that
\[\chi(t)=\lgd 1,&\quad 0\lsl t\lsl K_1,\\ 0,&\quad t\gsl K_2,\rgd\quad \text{and}\quad \chi'(t)<0,~K_1<t<K_2.\]

We define
\be H_K(z)=\chi(|z_{\sI}|)H(z)+(1-\chi(|z_{\sI}|))\rho|z_{\sI}|^{\mu},\ee
where
\be \rho\gsl \max_{K_1\lsl|z_{\sI}|\lsl K_2}\frac{|H(z)|}{|z_{\sI}|^\mu}.\lb{eq:defrho}\ee
The function $H_K$ satisfies (H0)-(H1) and
\[0<\mu H_K(z)\lsl z_\sI\cdot(H_K)_{z_\sI},~|z_\sI|\gsl r.\]
Then
\begin{align} &a_1|z_{\sI}|^\mu-a_2\lsl H_K(z)\lsl \rho|z_{\sI}|^\mu+a_2,\lb{eq:<HK<}\\
&\int^1_0\big(z_\sI\cdot (H_K)_{z_\sI}-H_K\big)\rd t\gsl\int^1_0\big((\mu-1)H-a_3\big)\rd t,\\
&H_K(z)<\min_{|z_\sI|\gsl K_2}H_K(z),\ |z_\sI|<K_2.\lb{eq:HK<HK.K2}
\end{align}
where
\be a_3=\max_{|z_\sI|\lsl r}\,(\mu-1)H(z)-\min_{|z_\sI|\lsl r}z_\sI\cdot H_{z_\sI}(z).\lb{eq:a3a4}\ee
Let
\be B_K(x)=\int^1_0 \big[TH_K(x(t)+t(0,v))+x(t)\cdot J(0,v)\big]\rd t,\lb{eq:defBK}\ee
and
\be \Phi_K(x)=A(x)-B_K(x),~\fa x\in\cE.\lb{eq:defPhiK}\ee
Certainly functionals $B_K$ and $\Phi_K$ are well-defined and of $C^1$ class,
and $\Phi_K$ satisfies (I1) and (I2). We next prove that $\Phi_K$ satisfies (I3) and (I4).


\bP\lb{prop:PhiK-I3} $\Phi_K$ satisfies (I3), i.e., (PS) condition.\eP

\pf  Let us consider a sequence
\[x^{(m)}=(w^{(m)},\tht^{(m)})\in\cH\ts T^k\]
such that
\be \Phi_K(w^{(m)},\tht^{(m)})\lsl c\quad \text{and}\quad \Phi'_K(w^{(m)},\tht^{(m)})\to 0
\quad\text{as}\quad m\ra\ift.\lb{eq:PS}\ee
Certainly $\big\{\tht^{(m)}\big\}_{m\in\N}$ has a convergent subsequence. We show that
$\big\{ w^{(m)}\big\}_{m\in\N}$ is bounded in $\cH$.

Let $w=w^{(m)},~\tht=\tht^{(m)}$, and decompose $w$ as
\be w=w^0_\sI+w^++w^-=w^0_{\sI}+\wtd w_{\sI}+\wtd w_{\sII}=w_\sI+\wtd w_\sII,\lb{eq:w.dec}\ee
where
\[w^\pm\in E^\pm,~w^0_{\sI}\in E^0_{\sI},~\wtd w_{\sI}\in E_{\sI},~\wtd w_{\sII}\in E_{\sII}.\]

By \eqref{eq:PS} and for large $m$ we have
\be \Phi_K(w,\tht)-\<\Phi'_K(w,\tht),w_{\sI}\>\lsl c+||w_{\sI}||.\lb{eq:Phi-Phi'wI}\ee

{\bf Case 1.} $k=n$.
\[w_{\sI}=w^0_\sI+\wtd w_\sI=(p,0)\in E^0_\sI\op E_\sI,~\wtd w_{\sII}=(0,q)\in E_\sII.\]
Set $\xi(t)=x(t)+t(0,v)$, we have
\begin{align} 
\Phi_K(x)&=\int^1_0 \big[p(t)\cdot\dot q(t)-TH_K(\xi(t))+p(t)\cdot v\big]\rd t,\lb{eq:Phi.k=n}\\
\Phi'_K(x)&=\big(\dot q+v-T(H_K)_p,-\dot p-T(H_K)_q\big),\nn\\
\<\Phi'_K(x),(p,0)\>&=\int^1_0 \big[p(t)\cdot \dot q(t)+p(t)\cdot v-p(t)\cdot T(H_K)_p(\xi(t))\big]\rd t,\nn\\
\Phi_K(x)-\<\Phi'_K(x),(p,0)\>&=T\int^1_0\big[p(t)\cdot (H_K)_p(\xi(t))-H_K(\xi(t))\big]\rd t\nn\\
&\gsl T\lt(\int^1_0 (\mu-1)H_K(\xi(t))\rd t-a_3\rt).\lb{eq:Phi-Phi'.k=n}
\end{align}
By \eqref{eq:Phi-Phi'wI} and \eqref{eq:Phi-Phi'.k=n}, we have
\[\bgd c+||w_{\sI}||&\gsl T\lt(\int^1_0(\mu-1)H_K(\xi(t))\rd t-a_3\rt)\\
&\gsl T\lt(\int^1_0(\mu-1)(a_1|w_{\sI}|^\mu-a_2)\rd t-a_3\rt).\egd\]
Then
\be||w_{\sI}||^\mu_\mu\lsl b_1(1+||w_{\sI}||)\lsl b_1(1+||w||).\ee
Here $||z||_\lm$ denotes the standard norm of $z\in L^\lm$. We have
\be |w^0_{\sI}|\lsl b_2(1+||w||^{1/\mu}).\lb{eq:p0<b2(1+w)}\ee

We next estimate $\wtd w_{\sI}+\wtd w_{\sII}=w^++w^-$ with $w^\pm\in E^\pm$. Note that
\be \<\Phi'_K(x),x^\pm\>=\int^1_0 \big[(x^+-x^-)\cdot x^\pm-TH'_K(\xi(t))\cdot x^\pm-J(0,v)\cdot x^\pm\big]\rd t.\lb{eq:Phi'xpm}\ee
By \eqref{eq:PS} we have
\[||w^\pm||^2\lsl\lt|\int^1_0TH'_K(\xi(t))\cdot w^\pm\rd t\rt|+||w^\pm||.\]
By definition of $H_K$ we have
\[|H'_K(z)|\lsl b_3(|z_\sI|^{\mu-1}+1).\]
Then
\[\lt|\int^1_0TH'_K(\xi(t))\cdot w^\pm\rd t\rt|\lsl b_3\int^1_0(|w_\sI|^{\mu-1}+1)|w^\pm|\rd t
\lsl b_3\lt(||w_\sI||^{\mu-1}_{\mu}+1\rt)||w^\pm||,\]
and
\be||w^\pm||\lsl b_4\lt(1+||w_\sI||^{\mu-1}_{\mu}\rt)
\lsl b_5\lt(1+||w||^{\frac{\mu-1}{\mu}}\rt).\lb{eq:wpm<=b4wI}\ee
By \eqref{eq:p0<b2(1+w)} we have
\[||w||\lsl||w^0_{\sI}||+||w^+||+||w^-||\lsl b_6\lt(1+||w||^{\frac{\mu-1}{\mu}}+||w||^{\frac{1}{\mu}}\rt).\]
Since $\mu>1$, we conclude that the sequence $\big\{w^{(m)}\big\}$ is bounded in $\cH$.

{\bf Case 2.} $k>n$. Let
\be\lgd
&x=(p_{\sI},p_{\sII},q_{\sI},q_{\sII}),\\
&p_{\sI}=(p_1,\cd,p_{2n-k}),~p_{\sII}=(p_{2n-k+1},\cd,p_n),\\
&q_{\sI}=(q_1,\cd,q_{2n-k}),~q_{\sII}=(q_{2n-k+1},\cd,q_n).\rgd\lb{eq:xI.xII.x.k>n}\ee
Then
\[w_\sI=(p_\sI,0,0,0),~w_\sII=(0,p_\sII,q_\sI,q_\sII).\]
Each vector $v\in\Z^k$ can be written as
\be v=(v'_{\sII},v_{\sI},v_{\sII}),~v_{\sI}\in\Z^{2n-k},~v_{\sII},v'_{\sII}\in\Z^{k-n}.\ee
We have
\begin{align}
&\Phi_K(x)=\int^1_0 \big(p_{\sI}\cdot\dot q_{\sI}+p_{\sII}\cdot \dot q_{\sII}
-TH_K(\xi(t))+p_{\sI}\cdot v_{\sI}+p_{\sII}\cdot v_{\sII}-q_{\sII}\cdot v'_{\sII}\big)\rd t,\\
&\Phi'_K(x)=\big(\dot q_{\sI}+v_{\sI}-T(H_K)_{p_{\sI}},\,\dot q_{\sII}+v_{\sII}-T(H_K)_{p_{\sII}},\nn\\
&\hs{1.8cm}-\dot p_{\sI}-T(H_K)_{q_{\sI}},\,-\dot p_{\sII}-v'_{\sII}-T(H_K)_{q_{\sII}}\big),\lb{eq:Phi'}\\
&\<\Phi'_K(x),(p_{\sI},0)\>=\int^1_0 \big(p_{\sI}\cdot\dot q_{\sI}+p_{\sI}\cdot v_{\sI}-p_{\sI}\cdot T(H_K)_{p_{\sI}}(\xi(t))\big)\rd t,\nn\\
&\Phi(x)-\<\Phi'(x),(p_{\sI},0)\>\nn\\
&=\int^1_0 \big(p_{\sI}\cdot T(H_K)_{p_{\sI}}(\xi(t))-TH_K(\xi(t))
+p_{\sII}\cdot(\dot q_{\sII}+v'_{\sII})-q_{\sII}\cdot v'_{\sII}\big)\rd t.\lb{eq:Phi-Phi'wI.2}\end{align}
By the definition of $H_K$, we have
\be|(H_K)_{z_{\sII}}|\lsl b_7.\lb{eq:HKzII<b7}\ee
By \eqref{eq:Phi'}, \eqref{eq:PS} and for large $m$ we have
\be||\dot p_{\sII}+v'_{\sII}+T(H_K)_{q_{\sII}}(\xi(t))||_2\lsl1,~||\dot q_{\sII}+v_{\sII}-(TH_K)_{p_{\sII}}(\xi(t))||_2\lsl 1.\lb{eq:2.43}\ee
 Since $w_{\sII}=(0,p_{\sII},q_{\sI},q_{\sII})\in E_{\sII}$, by Wirtinger's inequality we have
\[||p_{\sII}||_2\lsl\frac{1}{2\pi}||\dot p_{\sII}||_2,~||q_{\sII}||_{L^2}\lsl\frac{1}{2\pi}||\dot q_{\sII}||_2.\]
By \eqref{eq:HKzII<b7} and \eqref{eq:2.43} We have
\be\lt|\int^1_0 \big(p_{\sII}\cdot(\dot q_{\sII}+v'_{\sII})-q_{\sII}\cdot v'_{\sII}\big)\rd t\rt|\lsl b_8.\lb{eq:2.44}\ee
By \eqref{eq:Phi-Phi'wI.2}, \eqref{eq:Phi-Phi'wI}, (H2) and the above inequality, we have
\[\bgd c+||w_{\sI}||&\gsl T\lt(\int^1_0(\mu-1)H_K(\xi(t))\rd t-a_3\rt)\\
&\gsl T\lt(\int^1_0(\mu-1)(a_1|w_{\sI}|^\mu-a_2)\rd t-a_3\rt).\egd\]
The rest part of proof is similar with Case 1.~\QED\\

\bP\lb{prop:PhiK-I4} $\Phi_K$ satisfies (I4).\eP

\pf 

{\bf Case 1.} $k=n$. By \eqref{eq:XY:k=n}, each $x=(w,\tht)\in Y\ts T^n$ satisfies $w=(0,q)$, then
\[\Phi_K\bigg|_{Y\ts T^k}(x)=\int^1_0 -TH_K(0,q(t)+\tht+tv)\rd t\gsl\bb.\]
Each $x=(w,\tht)\in X\ts T^n$ can be written as
\[w=w^-+w^0,\ \ w^-=(p^-,q^-),~w^0=(p^0,0).\]
\begin{align}&\Phi_K\big|_{X\ts T^k}(x)=-||w^-||^2+p^0\cdot v-\int^1_0 TH_K(p^-+p^0,q^-+tv+\tht)\rd t\nn\\
&\lsl -||w^-||^2+|w^0|\cdot|v|-\int^1_0 T(a_1|p^-+p^0|^\mu-a_2)\rd t.\lb{eq:PhiX<.1}\end{align}

{\bf Subcase 1.1.} $1<\mu<2$. Let $f(t)=\dfrac{(1+t)^\mu}{1+t^\mu},~t\gsl0$. We have
\[\bgd f'(t)&=\frac{\mu(1+t)^{\mu-1}(1+t^\mu)-(1+t)^\mu\mu t^{\mu-1}}{(1+t^\mu)^2}\\
&=\frac{\mu(1+t)^{\mu-1}}{(1+t^\mu)^2}(1-t^{\mu-1})>0~~\text{iff}~~t<1.\egd\]
Then $\max\limits_{t\gsl0}f(t)=f(1)=2^{\mu-1}$ and
\[|u_2|^\mu=|u_1+u_2-u_1|^\mu\lsl(|u_1+u_2|+|u_1|)^\mu\lsl2^{\mu-1}(|u_1+u_2|^\mu+|u_1|^\mu),~\fa u_1,u_2\in\R^n.\]
Choose $u_1=p^-$ and $u_2=p^0$, we have
\[|w^0|=|p^0|^\mu\lsl2^{\mu-1}(|p^-+p^0|^\mu+|p^-|^\mu)\lsl 2^{\mu-1}(|p^-+p^0|^\mu+|w^-|^\mu).\]
By \eqref{eq:PhiX<.1} we have
\begin{align}
&\Phi_K|_{X\ts T^k}(x)\nn\\
&\lsl-||w^-||^2+|w^0|\cdot|v|-Ta_1\int^1_0\lt(2^{1-\mu}|w^0|^\mu-|w^-|^\mu\rt)\rd t+Ta_2\nn\\
&=-||w^-||^2+Ta_1||w^-||^\mu_\mu+|w^0|\cdot|v|-2^{1-\mu}Ta_1|w^0|^\mu+Ta_2\nn\\
&\lsl-||w^-||^2+Ta_1||w^-||^\mu+|w^0|\cdot|v|-2^{1-\mu}Ta_1|w^0|^\mu+Ta_2.\lb{eq:PhiX<.1.1}\end{align}
For each $w=w^-+w^0\in Q$, we have $||w^-||^2+|w^0|^2=R^2$. Then
\[||w^-||\gsl\frac{1}{\sqrt{2}}R~~\text{or}~~|w^0|\gsl\frac{1}{\sqrt{2}}R.\]
By \eqref{eq:PhiX<.1.1}, since $1<\mu<2$, for large $R$, we have
\be\Phi_K|_{\pt Q\ts T^k}(x)\lsl\bb-1.\lb{eq:cAKQ}\ee

{\bf Subcase 1.2.} $\mu\gsl2$. Similarly to computations in subcase 1.1, we have
\[\bgd\Phi_K|_{X\ts T^k}(x)&\lsl-||w^-||^2+|w^0|\cdot|v|-T\lt(a_1\int^1_0|p^-+p^0|^\mu\rd t-a_2\rt)\\
&\lsl-||w^-||^2+|w^0|\cdot|v|-T\lt(a_1\lt(\int^1_0|p^-+p^0|^2\rd t\rt)^{\frac{\mu}{2}}-a_2\rt)\\
&=-||w^-||^2+|w^0|\cdot|v|-T\lt(a_1\lt(\int^1_0|p^-|^2+|p^0|^2\rd t\rt)^{\frac{\mu}{2}}-a_2\rt)\\
&\lsl|w^0|\cdot|v|-Ta_1|w^0|^\mu+Ta_2.\egd\]
Then
\[\Phi_K|_{\pt Q\ts T^k}(x)\lsl|w^0|\cdot|v|-Ta_1|w^0|^\mu+Ta_2\lsl\bb-1\]
for sufficiently large $R$.

{\bf Case 2.} $k>n$. By \eqref{eq:XY:k>n},
\[X=\lt(E^-\cap X_1\rt)\op E^0_{\sI},~Y=E_{\sII}.\]
Similarly to Case 1,  we have
\[\Phi_K\bigg|_{Y\ts T^k}(x)\gsl\bb,~\fa x\in Y\ts T^k,\]
and
\[\Phi_K|_{\pt Q\ts T^k}(x)\lsl\bb-1,~\fa x\in \pt Q\ts T^k.~\QED\]

By Propositions \rf{prop:PhiK-I3}, \rf{prop:PhiK-I4} and \rf{prop:sptt.fp},
the functional $\Phi_K$ possess at least $k$ distinct critical points.
Rest proofs for \rT{thm:rot.fp.k=n} and \rT{thm:rot.fp.k>n} are preliminary estimates
for these critical points.


\ssct{Proof of \rT{thm:rot.fp.k=n}}

Let

\be \ga=\lgd T^{\frac{2}{2-\mu}}\cdot\lt(\frac{\mu a_1}{2}\rt)^{2/(2-\mu)}+T^{-\frac{1}{\mu-1}}\cdot2\lt(\frac{|v|^\mu}{\mu a_1}\rt)^{1/(\mu-1)}\lt(1-\frac{1}{\mu}\rt)+Ta_2,&\quad 1<\mu<2,\\
T^{-\frac{1}{\mu-1}}\cdot\lt(\frac{|v|^\mu}{\mu a_1}\rt)^{1/(\mu-1)}\lt(1-\frac{1}{\mu}\rt)+Ta_2,&\quad \mu\gsl2.\rgd\ee

\bP If $x(t)=\big(p(t),q(t)\big)$ is a critical point of $\Phi_K$ with $\Phi_K(x)\lsl \ga$, then
\[|p(t)|\lsl K,~H_K(x(t)+t(0,v))=H(x(t)+t(0,v)),~\fa t\in\R,\]
and $x$ is a solution of (\rf{eq:fpH.x}).\eP


\pf By \eqref{eq:Phi-Phi'.k=n},
\[\bgd \ga&\gsl\Phi(x)=\Phi(x)-\<\Phi'(x),(p,0)\>\\
&\gsl T\lt(\int^1_0 (\mu-1)H_K(p,q+tv)\rd t-a_3\rt)\\
&=T\big((\mu-1)H_K(p,q+tv)-a_3\big)\gsl T\lt((\mu-1)a_1|p(t)|^\mu-a_2-a_3\rt).\egd\]
We have $|p(t)|\lsl K$.~\QED

\ssct{Proof of \rT{thm:rot.fp.k>n}}

\bP If $x(t)=\big(x_{\sI}(t),x_{\sII}(t)\big)$ is a critical point of $\Phi_K$ with $\Phi_K(x)\lsl \ga$, then
\[|x_{\sI}(t)|\lsl K_1,~H_K(x(t)+t(0,v))=H(x(t)+t(0,v)),~\fa t\in\R,\]
and $x$ is a solution of (\rf{eq:fpH.x}).\eP

\pf Set
\[\xi(t)=x(t)+t(0,v).\]
According to \eqref{eq:xI.xII.x.k>n}, we have
\begin{align} \ga\gsl \Phi_K(x)&=\int^1_0 \big[p_{\sI}(t)\cdot\dot q_{\sI}(t)+p_{\sII}(t)\cdot \dot q_{\sII}(t)-TH_K(\xi(t))+p_{\sI}(t)\cdot v_{\sI}
+p_{\sII}(t)\cdot 0-q_{\sII}(t)\cdot 0\big]\rd t\nn\\
&=\int^1_0 \big[p_{\sI}(t)(\dot q_{\sI}(t)+v_{\sI})+p_{\sII}(t)\cdot\dot q_{\sII}(t)-TH_K(\xi(t))\big]\rd t\nn\\
&=\int^1_0 \big[p_{\sI}(t)\cdot (TH_K)_{p_{\sI}}(\xi(t))-TH_K(\xi(t))+\wtd p_{\sII}(t)\cdot\dot q_{\sII}(t)\big]\rd t\nn\\
&\gsl\int^1_0\big[T\big((\mu-1)H_K(\xi(t))-a_3\big)+\wtd p_{\sII}(t)\cdot(TH_K)_{p_{\sII}}(\xi(t))\big]\rd t.
\end{align}

For arbitrary $t,t'\in[0,1]$, we have
\[\wtd p_{\sII}(t)=\int^1_0\lt(\wtd p_{\sII}(t)-\wtd p_{\sII}(t')\rt)\rd t'
=\int^1_0\lt(\int^t_{t'}\dot p_{\sII}(t'')\rd t''\rt)\rd t'.\]
Then
\[\bgd &|\wtd p_{\sII}(t)|=\lt|\int^1_0\lt(\int^t_{t'}\dot p_{\sII}(t'')\rd t''\rt)\rd t'\rt|
\lsl\int^1_0\lt(\int^t_{t'}\lt|\dot p_{\sII}(t'')\rt|\rd t''\rt)\rd t'\\
&\lsl\int^1_0|\dot p_{\sII}(t'')|\rd t''=\int^1_0|\dot p_{\sII}(t)|\rd t
=\int^1_0|(TH_K)_{q_{\sII}}(\xi(t))|\rd t.\egd\]
We have
\[\bgd &\lt|\int^1_0\wtd p_{\sII}\cdot(TH_K)_{p_{\sII}}(\xi(t))\rd t\rt|\lsl \int^1_0|\wtd p_{\sII}|\,|(TH_K)_{p_{\sII}}(\xi(t))|\rd t\\
&\lsl\int^1_0|(TH_K)_{q_{\sII}}(\xi(t))|\rd t\int^1_0|(TH_K)_{p_{\sII}}(\xi(t))|\rd t\lsl\lt(\int^1_0|T(H_K)_{z_{\sII}}(\xi(t))|\rd t\rt)^2.\egd\]
By definition of $H_K$ and (H3), we have
\[|(H_K)_{z_{\sII}}|=|\chi(|p_{\sI}|)H_{z_{\sII}}|\lsl a|p_{\sI}|^s+b.\]
Then
\[\lt|\int^1_0\wtd p_{\sII}\cdot(TH_K)_{p_{\sII}}(\xi(t))\rd t\rt|\lsl T^2\lt(\int^1_0(a|p_{\sI}|^s+b)\rd t\rt)^2
\lsl T^2\cdot 2\lt(a^2||p_{\sI}||^{2s}_s+b^2\rt)\rd t.\]


By \eqref{eq:<HK<}, we have
\begin{align}
\ga&\gsl \int^1_0\big[T\big((\mu-1)H_K(\xi(t))-a_3\big)\big]\rd t-T^22\lt(a^2||p_{\sI}||^{2s}_s+b^2\rt)\lb{eq:ga>HK.mu}\\
&\gsl\int^1_0 T\big((\mu-1)(a_1|p_{\sI}|^\mu-a_2)-a_3\big)\rd t-T^22\lt(a^2||p_{\sI}||^{2s}_s+b^2\rt)\nn\\
&=T(\mu-1)a_1||p_{\sI}||^\mu_\mu-T^22a^2||p_{\sI}||^{2s}_s-T(\mu-1)a_2-Ta_3-T^22b^2\nn\\
&\gsl T(\mu-1)a_1||p_{\sI}||^{\mu}_s-T^22a^2||p_{\sI}||^{2s}_s-T(\mu-1)a_2-Ta_3-T^22b^2.\lb{eq:C>p.mu-p.2s}
\end{align}

Let
\be C\xeq{\text{def}}\ga+T(\mu-1)a_2+Ta_3+T^22b^2.\lb{eq:C.1}\ee
Then \eqref{eq:C>p.mu-p.2s} is
\be C\gsl T(\mu-1)a_1||p_{\sI}||^\mu_\mu-T^22a^2||p_{\sI}||^{2s}_s.\lb{eq:C>lm.mu-lm.2s}\ee
We claim that there exists a constant $R_1>0$ such that
\be ||p_\sI||_s\lsl R_1.\ee
By \eqref{eq:ga>HK.mu}, we have
\[\bgd \ga+Ta_3+T^22b^2+T^22a^2R^{2s}_1&\gsl T(\mu-1)\int^1_0H_K(\xi(t))\rd t\ev T(\mu-1)H_K(\xi(t))\\
&\gsl T(\mu-1)(a_1|p_\sI(t)|^\mu-a_2),~\fa t.\\
|p_\sI(t)|&\lsl \lt(\frac{C+T^22a^2R^{2s}_1}{T(\mu-1)a_1}\rt)^{1/\mu}=K_1.\egd\]

{\bf Case 1.} $2s<\mu$. By \eqref{eq:C>p.mu-p.2s}, we have $||p_\sI||_s\lsl R_1$.\vs{.5em}

{\bf Case 2.} $2s=\mu$. For $T\in\lt(0,\frac{(\mu-1)a_1}{2a^2}\rt)$ and
\[R_1=\lt(\frac{\ga+T(\mu-1)a_2+Ta_3+T^22b^2}{T(\mu-1)a_1-T^22a^2}\rt)^\frac{1}{\mu},\]
we have $||p_\sI||_s\lsl R_1$.\vs{.5em}

{\bf Case 3.} $\mu<2s<2\mu-1$.\vs{.5em}

{\bf Subcase 3.1.} $1<\mu<2$. Let $\lm=||p_\sI||_s$ and
\begin{align} C&\xeq{\text{def}}\ga+T(\mu-1)a_2+Ta_3+T^22b^2\nn\\
&=T^{\frac{2}{2-\mu}}\cdot\lt(\frac{\mu a_1}{2}\rt)^{\frac{2}{2-\mu}}
+T^{-\frac{1}{\mu-1}}\cdot 2\lt(\frac{|v|^\mu}{\mu a_1}\rt)^{\frac{1}{\mu-1}}\lt(1-\frac{1}{\mu}\rt)\nn\\
&\qquad +T\mu a_2+Ta_3+T^22b^2.\lb{eq:C.1}\end{align}

Set
\be \vf(\lm)=A\lm^\aa-B\lm^\bb,~\aa<\bb.\lb{eq:phi.1}\ee
Then
\[\bgd
\vf'(\lm)&=A\aa\lm^{\aa-1}-B\bb\lm^{\bb-1}=0~\mRa{~~}~\lm_0=\lt(\frac{A\aa}{B\bb}\rt)^{1/(\bb-\aa)}.\\
\max\vf(\lm)&=\vf(\lm_0)=B\lm^\bb_0\lt(\frac{\bb}{\aa}-1\rt)=B\lt(\frac{A\aa}{B\bb}\rt)^{\frac{\bb}{\bb-\aa}}\lt(\frac{\bb}{\aa}-1\rt).\egd\]
Choose
\[A=T(\mu-1)a_1,~B=T^22a^2,~\lm_0=\lt(\frac{T(\mu-1)a_1\mu}{T^22a^22s}\rt)^{\frac{1}{2s-\mu}},\]
we have
\begin{align}
\max\vf(\lm)&=T(\mu-1)a_1\lm^\mu_0-T^22a^2\lm^{2s}_0\nn\\
&=T^22a^2\lt(\frac{T(\mu-1)a_1\mu}{T^22a^22s}\rt)^{\frac{2s}{2s-\mu}}\lt(\frac{2s}{\mu}-1\rt)\nn\\
&=T^2\lt(\frac{T}{T^2}\rt)^{\frac{2s}{2s-\mu}}\cdot 2a^2\lt(\frac{(\mu-1)a_1\mu}{2a^22s}\rt)^{\frac{2s}{2s-\mu}}\lt(\frac{2s}{\mu}-1\rt)\nn\\
&=T^{\frac{2s-2\mu}{2s-\mu}}\cdot D_0,\lb{eq:phi.max}\end{align}
where
\begin{align} D_0&=2a^2\lt(\frac{(\mu-1)a_1\mu}{2a^22s}\rt)^{\frac{2s}{2s-\mu}}\lt(\frac{2s}{\mu}-1\rt)\nn\\
&=(2a^2)^{\frac{-\mu}{2s-\mu}}\lt(\frac{(\mu-1)a_1\mu}{2s}\rt)^{\frac{2s}{2s-\mu}}\lt(\frac{2s}{\mu}-1\rt).\end{align}
Let
\be D_1=\lt(\frac{\mu a_1}{2}\rt)^{\frac{2}{2-\mu}},~D_2=2\lt(\frac{|v|^\mu}{\mu a_1}\rt)^{\frac{1}{\mu-1}}\lt(1-\frac{1}{\mu}\rt),~D_3=\mu a_2+a_3,~D_4=2b^2.\ee

There exists $\dl>0$ such that for $T\in(0,\dl)$,
\be T^{\frac{2}{2-\mu}}\cdot D_1+T^{-\frac{1}{\mu-1}}\cdot D_2+T\cdot D_3+T^2\cdot D_4<T^{\frac{2s-2\mu}{2s-\mu}}\cdot D_5,\lb{eq:T<dl}\ee
i.e.,
\[T^{\frac{2}{2-\mu}+\frac{1}{\mu-1}}\cdot D_1+D_2+T^{1+\frac{1}{\mu-1}}\cdot D_3+T^{2+\frac{1}{\mu-1}}\cdot D_4<T^{\frac{\mu(2s-2\mu+1)}{(2s-\mu)(\mu-1)}}\cdot D_5.\]

In fact, since $\mu<2s<2\mu-1$,
\[\lls_{T\ra0^+}\text{Left}=D_2,~\lls_{T\ra0^+}\text{Right}=+\ift.\]
By \eqref{eq:C.1}, \eqref{eq:phi.1}, \eqref{eq:phi.max} and \eqref{eq:T<dl}
\begin{align} C&<T(\mu-1)a_1\lm^\mu_0-T^22a^2\lm^{2s}_0\nn\\
&=\max_{\lm\in\R}\lt(T(\mu-1)a_1\lm^\mu-T^22a^2\lm^{2s}\rt).\lb{eq:C<phi.max}
\end{align}
Then the inequality \eqref{eq:C>lm.mu-lm.2s} implies that
\be ||p_\sI||_s\lsl R_1~\text{or}~||p_\sI||_s\gsl R_2,~R_1<\lm_0<R_2.\lb{eq:R1<lm0<R2}\ee
By \eqref{eq:HK<HK.K2}, we have $||p_\sI||_s\lsl R_1$.~\QED

\bR In subcase 3.1, we must show that
\[K_1=\lt(\frac{C+T^22a^2R^{2s}_1}{T(\mu-1)a_1}\rt)^{1/\mu}<K_2=R_2.\]
In fact, by \eqref{eq:C<phi.max} and \eqref{eq:R1<lm0<R2}, we have
\[\bgd &R_1<\lm_0<R_2,\\
&C+T^22a^2\lm^{2s}_0<T(\mu-1)a_1\lm^\mu_0,\\
&K_1<\lt(\frac{C+T^22a^2\lm^{2s}_0}{T(\mu-1)a_1}\rt)^{1/\mu}<\lm_0<R_2=K_2.\egd\]
\eR



\end{document}